\documentclass[11pt]{amsart}

\usepackage{amsfonts,amssymb,amsopn,amscd,amsmath,amsthm}

\usepackage{graphicx}
\usepackage{color}

 \usepackage{mathrsfs}
 \numberwithin{equation}{section}
\theoremstyle{plain}
\newtheorem{thm}{Theorem}[section]
\newtheorem{lem}[thm]{Lemma}
\newtheorem{prop}[thm]{Proposition}

\newtheorem{rem}[thm]{Remark}
\theoremstyle{definition}

\newtheorem{ass}[thm]{Assumption}
\newtheorem{ex}{Example}[section]

\def\x{\mathbf{x}}

\def\X{\mathbf{X}}
\def\Y{\mathbf{Y}}
\def\Z{\mathbf{Z}}

\def\P{\mathbf{P}}
\def\A{\mathbf{A}}
\def\C{\mathbf{C}}

\def\R{\mathbb{R}}
\def\N{\mathbb{N}}

\def\z{\mathbf{z}}
\def\om{\mathbf{\Omega}}

\def\e{\mathbf{e}}
\def\x{\mathbf{x}}
\def\y{\mathbf{y}}
\def\h{\mathbf{h}}

\def\c{\mathbf{c}}

\def\V{\mathbf{V}}

\def\balpha{\boldsymbol{\alpha}}
\def\blambda{\boldsymbol{\lambda}}

\def\y{\mathbf{y}}

\def\0{\mathbf{0}}

\def\1{\mathbf{1}}

\begin{document}
\title[Max-entropy and computational geometry]
{Connecting Max-entropy with computational geometry, LP and SDP}
\author{Jean B. Lasserre}
\thanks{J.B. Lasserre is supported by the AI Interdisciplinary Institute ANITI  funding through the 
ANITI AI Cluster program under the Grant agreement ANR-23-IACL-0002, and also by the Marie-Sklodovska-Curie european doctoral network TENORS, under grant 101120296.
This research is also part of the programme DesCartes and is supported by the National Research Foundation, Prime Minister's Office, Singapore under its Campus for Research Excellence and Technological Enterprise (CREATE) programme.}

\address{LAAS-CNRS \& Toulouse School of Economics (TSE)}
\email{lasserre@laas.fr}
\begin{abstract}
 We consider the well-known max- (relative) entropy problem 
  $\Theta(\y)=\inf_{Q\ll P}D_{KL}(Q\Vert P)$ with Kullback-Leibler 
 divergence on a domain $\om\subset \R^d$,
 and with ``moment" constraints $\int\h\,dQ=\y$, $\y\in\R^m$.
 We show that when $m\leq d$, $\Theta$ is the Cram\'er transform of a function $v$
 that solves a simply related computational geometry problem. 
 Also, and remarkably,
  to the canonical LP:  $\min_{\x\geq0} \{\c^T\x: \A\,\x=\y\}$, with $\A\in\R^{m\times d}$,
 one may associate a max-entropy problem with a suitably chosen reference measure $P$ on $\R^d_+$ and 
 linear mapping $\h(\x)=\A\x$, such that its associated perspective function $\varepsilon\,\Theta(\y/\varepsilon)$
 is the optimal value of the log-barrier formulation (with parameter $\varepsilon$) of the dual LP (and so
  it converges to the LP optimal value as $\varepsilon\to 0$).
 An analogous result also holds for the canonical SDP: $\min_{X\succeq0} \{\langle C,X\rangle: 
 \mathcal{A}(X)=\y\,\}$.\\
 
 \end{abstract}
\maketitle

\section{Introduction}

The Max-(relative) entropy problem with Kullback-Leibler divergence criterion (let us 
call it KL-Maxent in the rest of the paper) is 
 important in statistics and in convex optimization on spaces of measures. 
Given a reference measure $P$ on $\om\subset \R^d$, the KL-divergence criterion
is used to select an optimal measure among all measures on $\om$ that satisfy $m$ generalized moment constraints.

In addition to being
an important problem in its own, KL-Maxent has a strong connection with the \emph{Large Deviation Principle} (LDP)
in probability. Indeed, under suitable conditions 
the optimal value function $\Theta$ associated with KL-Maxent is the \emph{Cram\'er rate function} in the 
LDP for  the empirical mean of a certain pushforward measure $\nu$. Equivalently, $\Theta$ is the \emph{Cram\'er transform} of the measure $\nu$; see for instance \cite{ben-tal,brown,Cramer,dembo,Donsker,gray,hollander,Yakov}. 
For instance in \cite{Yakov} the authors consider the particular case of  KL-Maxent on the mean (MEM), i.e. 
when the mean of the unknown distribution is imposed; they
describe (and prove some)  properties of the Cram\'er rate function, in particular in the context of exponential families. Also, depending on the chosen reference distribution of the model,
they also provide a quite informative list of Cram\'er rate functions in explicit form. 

In full generality,
solving KL-Maxent numerically is a computational challenge and for the interested reader, this issue is discussed 
in e.g. \cite{bach,jordan,Straszak}. However in \cite{Yakov} for the particular and favorable case MEM of KL-Maxent, efficient algorithms are proposed. For the interested reader, MEM has been used 
for solving inverse problems in Crystallography in \cite{paper-5,paper-6}, for  linear inverse problems 
(with approximate data-driven priors) in \cite{paper-7}, while the KL-divergence has been used as a regularizer in convex optimization for image processing (deblurring, denoising) in \cite{paper-2,paper-3}, and for a generalized inverse problem in computational Thermodynamics in \cite{paper-4}.

\vspace{.2cm}

While in the statistics and LDP literature
one speaks about the Cram\'er transform of \emph{measures} and not \emph{functions}, on the other hand the Cram\'er 
transform of functions is also an important tool of convex optimization; see e.g. the illuminating discussion in 
\cite[\S 9.4]{linearity}.\\

\noindent
{\bf Informal contribution.}
  
I.  In a first contribution of this  paper, we proceed further for KL-Maxent when $m\leq d$, and
provide a rigorous and explicit link between
KL-Maxent on the one side, and a related $(d-m)$-dimensional \emph{volume-type} computational geometry problem on the other side. 
Indeed when $m\leq d$,  the pushforward measure $\nu$ alluded to above, has a density with respect to 
(w.r.t.) Lebesgue measure on $\R^m$. Then by using geometric measure arguments (essentially Federer's coarea formula), 
one exhibits the density of $\nu$ as the value function $v$ of a simply related $(d-m)$-dimensional computational geometry problem. As a result, $\Theta$ now becomes the Cram\'er transform of that function $v$ which has a  
more familiar and intuitive flavor 
than the more abstract pushforward measure $\nu$. In particular, the partition-function 
$Z(\blambda)$ associated with  KL-Maxent is nothing less than $\mathcal{L}_v(-\blambda)$ ($=:\mathcal{L}_v^\vee(\blambda)$),
the Laplace transform of $v$ evaluated at $-\blambda$.
And indeed this view point is only valid when $m\leq d$, 
as $v$ does not  exist whenever $m>d$.  

Our motivation for considering the case $m\leq d$ comes from problems where
(i) $d$ is quite large and (ii) one  has access to only a few generalized moments $\y\in\R^m$. For instance, when the variables $\x\in\R^d$
are associated with nodes in a graph $G=(V,E)$ (e.g. as  in graphical models \cite{jordan}) 
with maximal cliques $\mathcal{C}_1,\ldots \mathcal{C}_p$, and the available information is 
only marginal through generalized moments associated with ``local" functions $h_{jk}:\R^{\#\mathcal{C}_j}\to\R$, $k=1,\ldots,r_j$, and $d\geq \sum_{j=1}^pr_j=:m$.

There are also several examples of cases $m\leq d$, especially when
$(m,\h)=(d,\mathrm{I}_d)$, referred to as \emph{Maximum Entropy on the Mean} (MEM), which has been
analyzed in detail in the recent work \cite{Yakov}. MEM was also considered in previous works 
\cite{paper-2,paper-3,paper-4,paper-5,paper-6,paper-7} for applications in different area
(e.g. in Crystallography and computational Thermodynamics). In particular, the authors in \cite{Yakov} provide a detailed discussion on sufficient conditions to ensure equivalence between $\Theta$ ($\kappa_p$ in \cite{Yakov}) and the  Cram\'er rate function (e.g., depending on whether the support of the reference distribution $P$ is countable or uncountable). They also list Cram\'er rate functions in closed form associated with several popular reference distributions \cite[Table 1, p. 457]{Yakov}, and provide a generic algorithm for MEM.

II. In a second contribution, one considers the canonical linear program (LP):
$\tau(\y)=\min_{\x\geq0}\{ \c^T\x: \A \x=\y\}$
(resp. the canonical semidefinite program (SDP):
$\tau(\y)=\min_{X\succeq0}\{ \langle C,X\rangle: \mathcal{A}(X)=\y\}$).
Then one may associate to  the LP (resp. the SDP) a KL-Maxent problem with an appropriate 
reference measure related to $\c$ (resp. $C$), and our previous general result specializes in a remarkable manner.
Indeed we show that the \emph{perspective function} $\varepsilon\,\Theta(\y/\varepsilon)$ 
(see e.g. \cite{perspective}) associated with
the optimal value function $\Theta$ of this KL-Maxent, is the same (up to a constant) as the optimal value of the log-barrier formulation with parameter $\varepsilon>0$, of the dual LP (resp. dual SDP). Hence
$\varepsilon\,\Theta(\y/\varepsilon)\to \tau(\y)$ as $\varepsilon\downarrow 0$. In addition,
$\Theta(\y)$ is obtained by solving the log-barrier formulation of the LP with parameter $\varepsilon=1$.
A  similar result also holds for the canonical SDP. To the best of the author's knowledge, this rigorous connection between the canonical LP (SDP) and a specific associated KL-Maxent problem is new.

Links between optimization and integration viewed
as formally the same ``operation" but in different algebra,
have been already discussed in several contexts; see e.g. \cite[\S 9.4]{linearity}, \cite{lasserre-book}, and the many references therein. 
In this paper, one provides another instance of this formal link (via Cram\'er transform)
between optimization and integration.

 \section{Setting, background \& contribution}
Let $\om\subset\R^d$ be open and with $m\leq d$, let
$\h:\om\to\R^m$ be continuously differentiable and Lipschitz. Let $j_\h(\x)=\sqrt{\mathrm{det}\,J_\h(\x)J_h(\x)^T}$ with $J_\h$ being the Jacobian associated with $\h$, and let 
$\mathcal{H}^{d-m}$ denote the $(d-m)$-dimensional Hausdorff measure. In the sequel one assumes that $j_\h>0$ for all $\x\in\om$.

Denote by $\mathscr{P}(\om)$ the set of probability measures on $\om$. With $\y\in\R^m$ and $P\in\mathscr{P}(\om)$ a given reference distribution  with support $\mathrm{supp}(P)=\overline{\om}$  and  density $p$ w.r.t. Lebesgue measure, consider the 
two computational problems :
\begin{eqnarray}
\label{pb-2}
v(\y)&:=&\int_{\om \cap \h^{-1}(\y)}\frac{p(\x)}{j_\h(\x)}\,d\mathcal{H}^{d-m}(\x)\\
\nonumber
\Theta(\y)&:=&\inf_{Q\ll P}\{\,\int_\om q(\x)\ln(\frac{q(\x)}{p(\x)})\,d\x\:: \int_\om \h(\x)\,dQ(\x)\,=\,\y\,\}\\
\label{pb-3}
&=&\inf_{Q\ll P}\{\,D_{KL}(Q\Vert P):\: \int_\om \h(\x)\,dQ(\x)\,=\,\y\,\}\,,
\end{eqnarray}
where the notation
$Q\ll P$ stands for ``$Q=q(\x)\,d\x$ is a probability measure on $\om$, absolutely continuous w.r.t. $P$". 
In \eqref{pb-2}
$\mathcal{H}^{d-m}$ is the  natural ``area" (or ``volume") measure on the fiber
$\h^{-1}(\y)$, well-defined as $\h$ is Lipschitz.

While Problem \eqref{pb-2} has a \emph{computational geometry}
flavor\footnote{Indeed, introducing the measure $V(d\x):=p/j_\h \mathcal{H}^{d-m}(d\x)$,
$v(\y)$ measures the ``size" of the set $\om\cap\h^{-1}(\y)=\{\x\in\om: \h(\x)=\y\}$, i.e.,
$v(\y)=V(\om\cap\h^{-1}(\y))$.}
  (e.g., when   $p\equiv j_\h$ in \eqref{pb-2}, $v$ is a $(d-m)$-dimensional \emph{volume}), 
 on the other 
hand Problem \eqref{pb-3} has a \emph{statistics} and \emph{information theory}  flavor as in \eqref{pb-3} one recognizes the 
\emph{maximum relative entropy} problem with Kullback-Leibler divergence
$D_{KL}(Q\Vert P)$ for a pair $(Q,P)$ of probability measures (let us call it the KL-Maxent problem).
Moreover, in problem \eqref{pb-3} the constraints (\emph{linear} in the unknown distribution $Q$), are
\emph{generalized} moment constraints (as $\h$ is not assumed to be polynomial) and so
\eqref{pb-3} addresses the Generalized Moment Problem by using the KL-divergence 
criterion to
select one among feasible solutions (if there is any feasible $Q\ll P$). 
Generalization to other \emph{$f$-divergences} have also been considered, see e.g. \cite{bach}, and
in \cite{ben-tal} max-entropy problems even more general  than \eqref{pb-3} are considered
by including additional (non linear but convex) entropy constraints. 
When $(m,\h)=(d,I_d)$ (i.e., the MEM problem) then trivially $v(\y)=p(\y)$ for all $\y\in\h(\om)$
as $\om\cap\h^{-1}\y)=\y$ whenever $\y\in\h(\om)$; indeed $\mathcal{H}^{0}$ is the counting measure, and $j_h(\x)=1$.

One would like to draw attention to an important difference between \eqref{pb-2} and \eqref{pb-3}.
Namely, while $v(\y)=0$ whenever
$\y\not\in \h(\om)$, on the other hand $\Theta(\y)$ is defined in 
the larger set $\overline{\mathcal{M}}$ (where $\mathcal{M}\:=\{\,\int \h\,dQ:\:Q\in\mathscr{P}(\om)\}$), hence in particular for $\y\in\mathrm{conv}(\h(\om))\setminus \h(\om)$.
However we will see that $v$ and $\Theta$ are still related in a precise mathematical sense.

\subsection*{Some background on KL-Maxent}
Recall that the \emph{Cram\'er transform} $\mathrm{Cramer}_f$ of a function $f$ is the Legendre-Fenchel transform $\mathscr{F}$
of the logarithm of its Laplace transform $\mathcal{L}_f$, that is:
\[\mathrm{Cramer}_f(\cdot)\,:=\,\mathscr{F}_{\ln \mathcal{L}_f}(\cdot)\,=\,
(\ln \mathcal{L}_f)^*(\cdot)\,,\]
(where $f^*$ is the usual notation for the Legendre-Fenchel transform of $f$)
and it plays an important role in convex optimization. 
See for instance the illuminating discussion in \cite[\S 9.4]{linearity} where 
seemingly different problems can be in fact seen as ``equivalent" problems in different algebra.
Depending on the context,
and for real $\blambda$, the Cram\'er transform either uses the Laplace transform $\mathcal{L}_f(\blambda)$ or the moment generating function (MGF) 
$\mathcal{L}^\vee_f(\blambda)\:(:=\mathcal{L}_f(-\blambda))$
(as e.g. in the statistics, information theory, and  LDP literature).

Associated with Problem \eqref{pb-3} is the 
\emph{log-partition} function 
$\ln Z(\blambda)
:=\ln\int\exp(\langle\blambda,\h\rangle)dP$)
an important tool in statistics and exponential families
\cite{bach,hollander,jordan}.  Under some assumptions on $\y$, $Z(\blambda)$ satisfies:
\begin{equation}
\label{infact-dual}
\mathscr{F}_{\ln Z}(\y)\,=\,(\ln Z)^*(\y)
\,:=\,\sup_{\blambda}\, \langle \blambda,\y\rangle-\ln Z(\blambda)\,=\,\Theta(\y)\,. 
\end{equation}
A typical condition to ensure that \eqref{infact-dual} holds 
is that $\y$ should belong to the relative interior of a certain convex set;
see e.g.  \cite{lewis,csiszar,dembo,jordan}. In the setting considered in this paper,
it is the convex set $\mathcal{M}$ introduced earlier, and
when $\h(\x)=(\x^{\balpha})_{\balpha\in\N^d;\vert\balpha\vert\leq n}$,
$\mathcal{M}$ is a slice of the so-called (truncated) \emph{moment cone}\footnote{With a closed set $K\subset\R^d$,
the moment cone associated with the degree-$n$ truncated $K$-moment problem is
the convex cone $\{(\int\,x_1^{\alpha_1}\cdots x_d^{\alpha_d}d\phi)_{\sum_i\alpha_i\leq n}: \phi\in\mathscr{M}(K)\}$
where $\mathscr{M}(K)$ is the convex cone of (positive) measures on $K$.}
in the literature on the moment problem (e.g. \cite{Schmudgen} and
\cite{didier} on a recent numerical scheme to test membership).

The maximization problem in \eqref{infact-dual} (which computes the Legendre-Fenchel transform of $\ln Z$)
is the \emph{dual} problem of \eqref{pb-3},
and \eqref{infact-dual} states that there is no duality gap between \eqref{pb-3} and its dual \eqref{infact-dual}.
When the ``$\inf$" in \eqref{pb-3} is attained at $Q^*$ and $\blambda^*$ is a maximizer 
in \eqref{infact-dual},  then 
\begin{equation}
\label{intro-Z}
Z(\blambda^*)\,=\, 
\int_\om \e^{\langle\blambda^*,\h\rangle}p(\x)\,d\x\,=\,
\int_\om \e^{\langle\blambda^*,\h\rangle}\,dP(\x)\,,\end{equation}
and 
\[Q^*(d\x)\,=\,q^*(\x)\,P(d\x)\,=\frac{\exp(\langle\blambda^*,\h(\x)\rangle)\,p(\x)\,d\x}{Z(\blambda^*)}\,.\]
For a discussion on the historical background about the equivalence of \eqref{pb-3} and \eqref{infact-dual}, the interested reader is referred to e.g. \cite{lewis,brown,csiszar,dembo,Donsker,hollander}, the recent \cite[\S 3]{Yakov} when $\h(\x)=\x$, and the many references therein. 
So \eqref{intro-Z} can be rewritten as
\begin{equation}
\label{intro-Z-1}
Z(\blambda^*)
\,=\,\int_{\h(\om)} \e^{\langle\blambda^*,\y\rangle}\h\#P(d\y)\,,\end{equation}
which is the MGF of the pushforward $\h\#P$ of $P$ by  the mapping $\h$.

The function $\Theta(\y)$ in \eqref{infact-dual} 
is called the \emph{Cram\'er 
rate function} in LDP for the empirical mean 
of the pushforward $\h\#P$ (and also denoted $I(\y)$
{\color{blue}in} the LDP literature);
see e.g. \cite{Cramer,dembo,Donsker,hollander} and more recently \cite{Yakov}. 

\begin{rem}
\label{rem-1}
Observe that when $m>d$, the pushforward $\h\#P$ of $P$ has \emph{no} density 
w.r.t. Lebesgue measure on $\R^m$ because its support 
$\overline{\h(\mathrm{supp}(P))}\subset\overline{\h(\om)}$ has  zero Lebesgue measure on $\R^m$. So the Cram\'er rate function $\Theta$ is \emph{not} the Cram\'er transform of 
some function $v:\R^m\to\R$. In the statistics and LDP literature,  the Cram\'er transform is defined for
probability distributions, not for functions. When $m>d$ one may still define a function $v$
but $v$ vanishes a.e. on $\R^m$ and integration of $v$ should be understood as integration 
w.r.t. $\mathcal{H}^d$ (on a submanifold of $\R^m$). For a precise  meaning, the interested reader is referred to e.g. \cite[Theorem 3.3]{Negro}.
 \end{rem}
 \subsection*{Contribution} 
However, when $m\leq d$ then one shows that in fact the computational geometry problem
\eqref{pb-2} is closely related to the statistical KL-Maxent problem \eqref{pb-3} through the Cram\'er transform of the function $v$ in \eqref{pb-2}. More precisely:

(i) With $v$ defined in \eqref{pb-2}, 
we show that $v$ is the density of the push-forward $\h\#P$ of $P$ by the mapping $\h$, and therefore
\eqref{intro-Z}-\eqref{intro-Z-1} 
also read
\begin{equation}
\label{intro-L_v}
Z(\blambda^*)\,=\,\int_{\h(\om)}\e^{\langle\blambda^*,\y\rangle}v(\y)\,d\y\,=:\,\mathcal{L}_v(-\blambda^*)\:(=:\,\mathcal{L}^\vee_v(\blambda^*))\end{equation}
(as $v=0$ whenever $\y\not \in \h(\om)$). As a consequence, from  \eqref{infact-dual}
\[\Theta(\y)\,=\,(\ln Z)^*(\y)\,=\,(\ln \mathcal{L}^\vee_v)^*(\y)\,=:\,\mathrm{Cramer}_v(\y)\,.\]
In other words,  $\Theta$ which is the \emph{Cram\'er rate function} associated with the empirical mean of
$\h\#P$ in LDP,
is also the \emph{Cram\'er transform} of the function $v$ in \eqref{pb-2}.

So even though the respective domains $\h(\om)$ of $v$,
and $\overline{\mathrm{conv}(\h(\om))}$ of $\Theta$ differ, $v$ and $\Theta$ are still related in a strong 
mathematical sense. In fact, $v$ encodes $\Theta$ via its Cram\'er transform. 
Hence when $m\leq d$, arguments from geometric measure theory (essentially Federer's coarea formula) allow to identify the rate function $\Theta$ as the Cram\'er transform of $v$ in \eqref{pb-2}, i.e.,
of the density of $\h\#P$ w.r.t. Lebesgue measure in $\R^m$,
which is well-defined if $j_\h>0$ on $\om$ and $\h$ is Lipschitz.

To illustrate the link between \eqref{pb-2} and \eqref{pb-3}, let
the reference distribution on $\om$ be $P_0(d\x)=j_\h(\x)d\x/s$, with $s:=\int_\om j_\h d\x$
(assuming $s<+\infty$). Then 
$v(\y)=s^{-1}\mathcal{H}^{d-m}(\om\cap \h^{-1}(\y))$ is the $(d-m)$-dimensional volume of the fiber $\om\cap\h^{-1}(\y)$ scaled by $1/s$. Therefore when using $P_0$ as reference probability distribution,
$\Theta(\y)$ is the Cram\'er transform of the fiber-volume function $v/s$.

(ii) Next, introduce the canonical linear program (LP) 
\begin{equation}
\label{lp-primal}
\tau(\y)=\min_{\x\geq0}\{\,\langle\c,\x\rangle: \A\x=\y\,\}\,
\end{equation}
(with $\c\in\R^d,\,\A\in\R^{m\times d}$,  and hence $m\leq d$), and its dual
\begin{equation}
\label{lp-dual}
\max_{\blambda} \,\{\,\langle\,\y,\blambda\rangle:\:\A^T\blambda \leq\c\,\}\,=\,\tau^*(\y)\:(=\tau(\y))\,.
\end{equation}
Assume that $\A^T\blambda_0<\c$ for some $\blambda_0$ so that 
$\tau(\y)>-\infty$ and finite if $\{\x\geq0:\A\x=\y\}\neq\emptyset$.
So with no loss of generality, one may and will assume that $\c>0$  (otherwise change $\c$ to $\tilde{\c}:=\c-\A^T\blambda_0>0$)\footnote{For a vector $\c\in\R^d$,  $\c>0$ stands for $c_i>$ for all $i=1,\ldots,d$.}. 

Then one shows that with the choice of reference probability measure 
$P(d\x)=s\,\e^{-\langle\c,\x\rangle}d\x$ on $\om=\R^d_+$ (with $1/s:=\int_\om \e^{-\langle\c,\x\rangle}d\x$), and $\h(\x):=\A\x$, the \emph{perspective function} $\varepsilon\,\Theta(\y/\varepsilon)$ (with $\varepsilon>0$ fixed) associated with $\Theta$ in \eqref{pb-3}, is also the optimal value (up to a constant) of
\begin{equation}
\label{log-barrier}
\tau_\varepsilon(\y)\,=\,\max_{\blambda}\,\{\,\langle\,\y,\blambda\rangle+\varepsilon\,\sum_{j=1}^m \ln (\c-\A^T\blambda)_j\,,\quad\A^T\blambda\,<\c\,\}\,,
\end{equation}
the classical log-barrier formulation of \eqref{lp-dual} for interior points methods.
As is well-known, $\tau_\varepsilon(\y)\to\tau(\y)$ as $\varepsilon\downarrow 0$, and therefore
$\varepsilon\,\Theta(\y/\varepsilon)\to \tau(\y)$ as $\varepsilon\downarrow 0$;
see e.g. \cite{etienne,Monteiro}.

So essentially,
with the particular choice of reference measure $P_0$, solving the KL-Maxent \eqref{pb-3} reduces to solving the LP
\eqref{lp-dual} via the log-barrier algorithm with parameter $\varepsilon=1$. Moreover,
in this case the function $v$ in \eqref{pb-2} has a nice explicit formulation in terms of (i) the vertices 
of the convex polyhedron $\{\x\geq0: \A\x=\y\}$, and (ii) reduced costs in the simplex algorithm, thanks 
to a formula of Brion \& Vergnes \cite{brion}.

(iii) A similar and analogous result also holds for the canonical semidefinite program (SDP)
\begin{equation}
\label{intro-sdp-primal}
\rho(\y)\,=\,\min_{\X\succeq0}\,\{\,\langle \C,\X\rangle:\: \langle \A_j,\X\rangle\,=\,y_j\,,\quad 1\leq j\leq m\,\}\,,
\end{equation}
and its dual
\begin{equation}
\label{intro-sdp-dual}
\sup_{\blambda} \,\{\,\langle\blambda,\y\rangle:\:\C-\sum_{j=1}^m\lambda_j\,\A_j\,\succeq0\,\}\,=\,\rho^*(\y)\,,
\end{equation}
for real symmetric matrices $\C,\A_j\in\R^{d\times d}$, and vector $\y\in\R^m$ with
$m <s_d:=d(d+1)/2$. 
Assume that $\C\succ\sum_{j=1}^m\lambda_{0j}\,
\A_j$ for some $\blambda_0\in\R^m$, so that 
$\rho(\y)$ is bounded below and finite if 
$\{\,\X\succeq0:\langle \A_j,\X\rangle=y_j,\:j\leq m\}\neq\emptyset$. 
Hence with no loss of generality, one may and will assume that $\C\succ0$ (otherwise change 
$\C$ to $\tilde{\C}:=\C-\sum_{j=1}^m\lambda_{0j}\A_j\succ0$). 

Let $\om:=\mathcal{S}^d_+\subset\mathcal{S}^d$ be the convex cone of real symmetric $(d\times d)$ positive semidefinite matrices,
and $\h:\mathcal{S}^d\to\R^m$ a linear mapping with $m<s_d$, and 
$\h(\x)=(\langle \A_j,\X\rangle)_{1\leq j\leq m}$.
Then we show that with the choice of reference probability measure 
$P(d\x)=s\,\e^{-\langle \C,\X\rangle}dX$ on $\om$ (with $1/s=\int_\om \e^{-\langle \C,\X\rangle}d\X$), the \emph{perspective function} $\varepsilon\,\Theta(\y/\varepsilon)$ (with $\varepsilon>0$ fixed) associated with $\Theta$, is also the optimal value (up to a constant) of
\begin{equation}
\label{log-barrier-sdp}
\rho^*_\varepsilon(\y)\,=\,\sup_{\blambda}\,\{\,\langle\blambda,\y\rangle+\varepsilon\,\sum_{j=1}^m 
\ln \mathrm{det}(\C-\sum_{j=1}^m \lambda_j\,\A_j)\,,
\end{equation}
the classical log-barrier formulation of \eqref{intro-sdp-dual} for
interior points methods. As is well-known \cite{log-barrier-3,etienne,Monteiro}, $\rho^*_\varepsilon(\y)\to\rho^*(\y)$ as $\varepsilon\downarrow 0$,
and therefore
$\varepsilon\,\Theta(\y/\varepsilon)\to \rho^*(\y)$ as $\varepsilon\downarrow 0$. So essentially,
with particular choice of reference measure $P$, solving the KL-Maxent \eqref{pb-2} reduces to solving the 
SDP \eqref{intro-sdp-dual} (and hence \eqref{intro-sdp-primal} as well) via the log-barrier algorithm with parameter $\varepsilon=1$.

\section{Main result}
\label{sec:main}
\subsection{Preliminaries}
\label{Laplace}
Let $\om\subset\R^d$ be open, and with $m\leq d$, 
let $h_1,\ldots,h_m:\R^d\to\R$
be continuously differentiable and Lipschitz.
Then with
$\h=(h_1,\ldots,h_m)^T:\R^d\to\R^m$, 
let $J_\h$ denote its Jacobian and $j_\h:\R^d\to\R_+$,
\[\x\mapsto j_\h(\x)\,:=\,\sqrt{\mathrm{det}(J_\h(\x)\,J_\h(\x)^T)}\,,\quad\x\in\R^d\,.\]

\subsection*{Multivariate (bilateral) Laplace transform}
Given $v:\R^d\to\R$, the function $\mathcal{L}_v:\mathbb{C}^d\to\mathbb{C}$, where
\[\mathcal{L}_v(\blambda)\,:=\,\int_{\R^d_+}\e^{-\langle \blambda,\y\rangle}\,v(\y)\,d\y\,,\quad\blambda\in\mathcal{D}\,,\]
is the  (bilateral) Laplace transform of $v$, with domain $\mathcal{D}\subset \mathbb{C}^d$ where the integral in the right-hand-side is finite.

\subsection*{Legendre-Fenchel transform}
Given $v:\R^d\to\R\cup\{+\infty\}$, the function
$v^*:\mathbb{R}^d\to\mathbb{R}\cup\{+\infty\}$ with
\[v^*(\blambda)\,:=\,\sup_{\y\in \R^d}\langle\blambda,\y\rangle-v(\y)\,,\quad\blambda\in\mathcal{D}\,,\]
denote its Legendre-Fenchel transform (or Fenchel conjugate) with domain $\mathcal{D}\subset \R^d$.

\subsection*{Hausdorff measure}
In geometric measure theory, the Hausdorff measure generalizes the Lebesgue measure and
the associated notions of counting, length, area and volume. 
Technically, let $S$
be any subset of a metric space $(X,\rho)$,
$\delta>0$ a real number, $d\in\N$, and define
\[H_{\delta }^{d}(S)\,:=\,\inf \left\{\sum _{i=1}^{\infty }(\mathrm{diam} \,U_i)^{d}:S\subseteq \bigcup _{i=1}^{\infty }U_{i},\:\mathrm{diam}\,U_i\,<\,\delta \right\}\,,\]
where the infimum is over all countable covers of  $S$
by sets $U_i\subset X$
satisfying $\mathrm{diam}\,U_i<\delta$. Then 
\[H^d(S)\,:=\,\sup_{\delta>0}H^d_\delta(S)\,=\,\lim_{\delta\downarrow 0}H^d_\delta(S)\]
exists (possibly $+\infty$).
The restriction 
of $H^d(\cdot)$ to sets $S$ in the $\sigma$-field of Caratheodory-measurable sets  is a measure called the $d$-dimensional Hausdorff measure of $S$. So the $d$-dimensional 
Hausdorff measure of $\R^d$ is a rescaling of 
the usual normalization of Lebesgue measure
(normalized to give volume $1$ to the unit box $[0,1]^d$). A standard reference is \cite{federer}.

\subsection*{Coarea formula}

With $\om$, $\h$ as above,
recall Federer's coarea formula \cite[Theorem 3.1, p. 426]{federer} for Lipschitz functions $\h$. 
Let $\mathcal{H}^{d-m}$ be the $(d-m)$-dimensional Hausdorff measure, and 
let $g:\R^d\to [0,+\infty)$ be measurable. Then:
\begin{equation}
 \label{coarea-0}
 \int_{\R^d} g(\x)\,j_\h(\x)\,d\x\,=\,\int_{\R^m}\left(\int_{ \h^{-1}(\y)}g(\x)\,d\mathcal{H}^{d-m}(\x)\right)\,d\y\,,
 \end{equation}
  provided that the integral on the left is finite (e.g. if $\om$ is bounded).
 Hence with $\x\mapsto g(\x)\1_\om(\x)$,
 \begin{equation}
 \label{coarea-1}
 \int_\om g(\x)\,j_\h(\x)\,d\x\,=\,\int_{\R^m}\left(\int_{\om\cap \h^{-1}(\y)}g(\x)\,d\mathcal{H}^{d-m}(\x)\right)\,d\y\,.
 \end{equation}
 In particular,
 \begin{equation}
 \label{coarea-2}
 \int_\om g(\h(\x))\,j_\h(\x)\,d\x\,=\,
 \int_{\R^m} g(\y)\,\mathcal{H}^{d-m}(\om\cap \h^{-1}(\y))
 \,d\y\,.
 \end{equation}
 \begin{ass}
 \label{ass-1}
 Throughout the rest of the paper:
 
 (i) $m\leq d$ and $j_\h$ is strictly positive
 on $\om$. 
 
 (ii) $p:\om\to\R_+$ is measurable and nonnegative.
 \end{ass}
Under Assumption \ref{ass-1},  letting $\x\mapsto g(\h(\x))\frac{p(\x)}{j_\h(\x)}$ in \eqref{coarea-1} yields:
  \begin{equation}
 \label{coarea-3}
 \int_\om g(\h(\x))\,p(\x)\,d\x\,=\,\int_{\R^m}g(\y)\,\left(\int_{\om\cap \h^{-1}(\y)}\frac{p(\x)}{j_\h(\x)}\,d\mathcal{H}^{d-m}(\x)\right)\,d\y\,
 \end{equation}
 whenever the integral on the left is finite (e.g. if $\om$ is bounded).
 \subsection*{The function $v$}  We here assume that $P$ is a finite  measure on $\om$ 
 with density $p$ w.r.t. Lebesgue measure on $\R^d$, and so $\int_\om p(\x)d\x=P(\om)$. Hence from \eqref{coarea-3} with $g\equiv 1$ (the constant function
 equal to $1$), and as $j_\h>0$ on $\om$,
 \[P(\om)\,=\,\int_\om p(\x)\,d\x\,=\,
 \int_{\R^m}\left(\int_{\om\cap\h^{-1}(\y)} \frac{p(\x)}{j_\h(\x)}\,d\mathcal{H}^{d-m}(\x)\right)\,d\y\,.\]
Next, introduce the function $v:\h(\om)\to\R_+$:
\begin{equation}
\label{def-v}
\y\mapsto v(\y)\,:=\,\int_{\om\cap \h^{-1}(\y)}\frac{p(\x)}{j_\h(\x)}\,d\mathcal{H}^{d-m}(\x)\,,\end{equation}
which is well-defined (possibly infinite) for all $\y
\in\h(\om)$ as we have assumed that $j_\h >0$ on $\om$. Next,
$v(\y)=0$ whenever $\mathcal{H}^{d-m}(\om\cap \h^{-1}(\y))=0$ (e.g. whenever 
$\y\not\in \h(\om)$), and since the left-hand-side of \eqref{def-v}
is finite then $v$ is finite for almost all $\y\in\R^m$.

 In particular, if $\int_{\om}j_\h(\x)d\x<+\infty$, then  for the choice 
$P(d\x)\,:=\,j_\h(\x)d\x$, one obtains
\[v(\y)\,:=\,\mathcal{H}^{d-m}(\om\cap \h^{-1}(\y))\,,\quad\mbox{a.a. }\y\in\h(\om)\,,\]
which is the $(d-m)$-dimensional ``volume" of the fiber $\om\cap \h^{-1}(\y)$.

\begin{prop}
\label{prop1}
Let Assumption \ref{ass-1} hold, with $P$ a finite measure on $\om$ with  density $p$ w.r.t. Lebesgue measure.
The measure $V(d\y)=v(\y)d\y$ on $\R^m$ with $v$ as in \eqref{def-v}, is the pushforward $\h_\#P$ of $P$ by the mapping $\h$. In particular, if $\int_\om j_\h(\x)\,d\x<+\infty$, then for the reference measure $P(d\x):=j_\h(\x)d\x$,
\begin{equation}
\label{prop1-1}
V(d\y)\,=\,\mathcal{H}^{d-m}(\om\cap\h^{-1}(\y))\,d\y\,=\,``\mathrm{vol}_{d-m}"(\om\cap \h^{-1}(\y))\,d\y\,.\end{equation}
\end{prop}
\begin{proof}
 $P$ is a finite measure. Let $A\in\mathcal{B}(\R^m)$ be fixed, arbitrary.
 Then
 \[\int_{\R^m}\1_A(\y)\,d\h\#P(\y)\,=\,\int_{\om}\1_A(\h(\x))\,dP(\x)
 \,=\,\int_{\om}\1_A(\h(\x))\,p(\x)\,d\x\,.\]
 Therefore, 
 \begin{eqnarray*}
\h \#P(A)&=&\int_\om\1_A(\h(\x))\,p(\x)\,d\x\\
&=&\int_{A}(\underbrace{\int_{\om\cap \h^{-1}(\y)}\frac{p(\x)}{j_\h(\x)}d\mathcal{H}^{d-m}(\x)}_{=v(\y)})\,d\y\,,
\end{eqnarray*}
 where we have used \eqref{coarea-1} with $\x\mapsto g(\x):=\1_A(\h(\x))p(\x)/j_\h(\x)$.  
 As this is valid for all Borel sets $A\in\mathcal{B}(\R^m)$,  one obtains the desired result 
 $V(d\y)=v(\y)d\y$ with $v$ in \eqref{def-v} being the density of $V$ for almost all $\y\in\R^m$. 
 
 Finally, for the particular reference 
 measure  $P(d\x)=j_\h(\x)d\x$ (when $\int_\om j_\h(\x)d\x<+\infty$), 
 in \eqref{def-v} the density of $V$ reads $v(\y)=\mathcal{H}^{d-m}(\om\cap\h^{-1}(\y))$
  for almost all $\y\in\R^m$, which is \eqref{prop1-1}.
   \end{proof}
\begin{rem}
\label{positivity}
In Assumption \ref{ass-1}(i), one may relax the condition of strict positivity of $j_\h$ on $\om$. Let $S:=\{\x\in\om: j_\h(\x)=0\}$, and define
\[\tilde{p}(\x)\,:=\,\left\{\begin{array}{rl}
\frac{p(\x)}{j_\h(\x)} &\mbox{if $\x\in\om\setminus S$}\\
0&\mbox{if $\x\in S$.}\end{array}\right.\]
Then as $S$ is a set of Lebesgue measure zero,
\begin{eqnarray*}
\h\#P(A)&=&\int_\om 1_A(\h(\x))\,p(d\x)d\x\,=\,\int_\om 1_A(\h(\x))\,\tilde{p}(\x)\,j_\h(\x)\,d\x\\
&=&\int_{\R^m}1_A(\y)\,\left(\int_{\om\cap\h^{-1}(\y)}\tilde{p}(\x)\,d\mathcal{H}^{d-m}(\x)\right)\,d\y\\
&=&\int_{\R^m}1_A(\y)\,
\left(\int_{(\om\setminus S)\cap\h^{-1}(\y)}\frac{p(\x)}{j_\h(\x)}\,
d\mathcal{H}^{d-m}(\x)\right.\,d\y\\
&&+\left. \int_{S\cap\h^{-1}(\y)}\underbrace{\tilde{p}(\x)}_{=0}\,d\mathcal{H}^{d-m}(\x)\right)\,d\y\,,
\end{eqnarray*}
for all Borel sets $A\in\mathcal{B}(\R^m)$, and so for almost all  $\y\in\R^m$,
\[\y\mapsto v(\y)\,=\,\int_{(\om\setminus S)\cap\h^{-1}(\y)}\frac{p(\x)}{j_\h(\x)}\,
d\mathcal{H}^{d-m}(\x)\,\]
is the density of the pushforward $\h\# P$. The difference 
with the case $j_\h(\x)>0$ for all $\x\in\om$, is that the definition of $v$ is now slightly more complicated as
it is an integral on $(\om\cap \{\x: j_\h(\x)>0\})\cap\h^{-1}(\y)$ instead of the simpler  set
$\om\cap\h^{-1}(\y)$.
\end{rem}
\begin{lem} 
\label{lem1}
Let Assumption \ref{ass-1} hold and let $v:\h(\om)\to\R_+$ be as in \eqref{def-v}. 
Then with $\blambda\in\R^m$,
\begin{eqnarray}
\label{lem1-1}
\int_\om\e^{\langle\blambda,\h(\x)\rangle}\,p(\x)\,d\x&=&
\int_{\R^m} \e^{\langle\blambda,\y\rangle}\,v(\y)\,d\y\,=\,
\mathcal{L}_v(-\blambda)\,(=:\mathcal{L}^\vee_v(\blambda))
\end{eqnarray}
whenever the integral on the left is finite. That is, the integral on the left is the \emph{Laplace transform}\ $\mathcal{L}_v$ of $v$ (evaluated at $-\blambda\in\R^m$) or the MGF $\mathcal{L}^\vee_v$ of $v$, evaluated at $\blambda$.
\end{lem}
\begin{proof}
With  $\x\mapsto \e^{\langle\blambda,\h(\x)\rangle}p(\x)/j_\h(\x)$ in  the coarea formula \eqref{coarea-1}-\eqref{coarea-2}, one obtains the desired result
  \begin{eqnarray*}
 \int_\om \e^{\langle\blambda,\h(\x)\rangle}\,p(\x)d\x&=&
 \int_{\R^m} \e^{\langle\blambda,\y\rangle}
 \left(\int_{\om\cap \h^{-1}(\y)}\frac{p}{j_\h}d\mathcal{H}^{d-m}\right)\,d\y\\
 &=&\int_{\R^m} \e^{\langle\blambda,\y\rangle}\,v(\y)\,d\y\,=\,\mathcal{L}_v(-\blambda)\,=\,\mathcal{L}^\vee_v(\blambda)\,,
 \end{eqnarray*}
 whenever the integral on the left-hand-side is finite.
 \end{proof}

\subsection{KL-Maxent (relative) entropy problem}

Let $\om\subset\R^d$ be open. Let $p,h_1,\ldots,h_m:\R^d\to\R$ be as in \S \ref{Laplace} and let Assumption \ref{ass-1} hold.  Denote by $\mathscr{P}(\om)$ be the set of probability measures on $\om$ and
assume that $\int_\om p(\x)d\x=1$ (i.e., $P\in\mathscr{P}(\om)$) and $\mathrm{supp}(\P)=\overline{\om}$. 
Let $\mathscr{C}_P(\om):=\{\,Q\in\mathscr{P}(\om):\:Q\ll P\,\}$ and write $Q(d\x)=q(\x)\,d\x$ 
with $q\in L^1(\om)_+$.
In addition define
\begin{equation}
\label{moment cone}
\mathcal{M}\,:=\,
\{\,\int_{\om}\h(\x)\,\mu(d\x)\,:\: \mu\in\mathscr{P}(\om)\,\}\,,
\end{equation}
i.e., $\mathcal{M}$ is the convex set of all realizable $\y\in\R^m$ by some probability distribution on $\om$
(not necessarily absolutely continuous w.r.t. $P$).

\subsection*{Kullback-Leibler divergence}
Let $P(d\x)=p(\x)d\x$ be a given reference probability measure with support $\overline{\om}$.
The KL -- (relative) entropy $D_{KL}(Q\Vert P)$ defined by
\[D_\mathrm{KL}(Q\Vert P)\,:=\,\left\{\begin{array}{rl}
\displaystyle\int_\om q(\x)\,\ln(\frac{q(\x)}{p(\x)})\,d\x&\mbox{if }Q\,\in\,\mathscr{C}_P(\om)\,,\\
+\infty&\mbox{otherwise,}\end{array}\right.\]
measures how much $Q$ \emph{diverges} from $P$, and is  called the \emph{Kullback-Leibler (KL) divergence}.

\subsection*{KL-Maxent problem}
With $\y\in\R^m$, the KL-Maxent problem is the convex optimization problem with optimal value function:
\begin{eqnarray}
\label{def-entropy}
\Theta(\y)&:=&\inf_{Q\in\mathscr{C}_P(\om)}\{\,\int_\om q(\x)\,\ln(\frac{q(\x)}{p(\x)})\,d\x
: \:\displaystyle\int_\om \h(\x)\,q(\x)\,d\x\,=\,\y\,\}\\
\nonumber
&=&\inf_{Q\in\mathscr{C}_P(\om)} \{\,D_\mathrm{KL}(Q\Vert P):\: \:\displaystyle\int_\om\h(\x)\,q(\x)\,d\x\,=\,\y\,\}\,.
\end{eqnarray}
The function $Z:\R^m\to\R_+$ defined by
\begin{equation}
\label{D_z}
Z(\blambda)\,:=\,\int_\om \exp(\langle\blambda,\h(\x))\,dP(\x)\,,
\end{equation}
with domain:
\begin{equation}
\label{domain-D_z}
\mathcal{D}_Z\,:=\,\{\,\blambda\in\R^m: Z(\blambda)<+\infty\,\}\,,
\end{equation}
is called the \emph{partition function}, and plays an important role in statistics, information theory and exponential families.  Notice that
by introducing the pushforward measure $\h\#P$ of
$P$ by the mapping $\h$, 
\eqref{D_z}  also reads
\begin{equation}
\label{Z(lambda)-push}
\lambda\mapsto Z(\blambda)\,:=\,\int_{\h(\om)} \e^{\langle\blambda,\y\rangle}\,\h\#P(d\y)\,.
\end{equation}
As the support of $\h\#P$ is $\overline{\h(\om)}$, $Z(\blambda)$ is the 
Moment Generating Function (MGF) of the pushforward distribution $\h\#P$  on $\h(\om)$ (or 
the Laplace transform $\h\#P$, evaluated at $-\blambda$) .
\begin{ass}
\label{ass-2}
$P(d\x)=p(\x)d\x$ is a probability measure with support $\overline{\om}$.
The domain $\mathcal{D}_Z$ in \eqref{domain-D_z} is open
and \emph{linear independence $P$-a.e.} holds, that is, if $\x\mapsto \langle\blambda,\h(\x)\rangle$ is constant, $P$-a.e., then  $\blambda=\mathbf{0}$.
\end{ass}
When Assumption \ref{ass-2} is interpreted in  the context of exponential families investigated in \cite[\S 3.5]{jordan}, 
the family of densities $(\exp(\langle\blambda,\h\rangle-\ln Z(\blambda))_{\blambda}$ w.r.t. $P$, form a \emph{regular} exponential family,
$\h$ is called a sufficient statistic, and the linear independence states that the exponential family is \emph{minimal}.  Then 
$\ln Z$ (also called the cumulant function) is strictly convex on its domain $\mathcal{D}_Z$ \cite[Proposition 3.1(b)]{jordan}. In addition, 
the gradient mapping $\nabla \ln Z:\mathcal{D}_Z\to\mathcal{M}$ is \emph{one-to-one} \cite[Proposition 3.2]{jordan},
and is \emph{onto} $\mathrm{int}(\mathcal{M})$  \cite[Theorem 3.3]{jordan},

As emphasized in \cite{jordan}, for $\y\in\mathrm{int}(\mathcal{M})$, there are many
probability distributions $\mu$ on $\om$ (not necessarily
member of an exponential family) such that $\int \h\,d\mu=\y$,
but the one that maximizes the entropy (i.e. minimizes \eqref{def-entropy}) is necessarily a member of the exponential family $\exp(\langle\blambda,\h\rangle-\ln Z(\blambda))_{\blambda}$.

Indeed, under Assumption \ref{ass-2}, 
\eqref{def-entropy} has a unique
optimal solution $q^*$, of the form
\begin{equation}
\label{def-optimal-q*}
\x\mapsto q^*(\x)\,:=\,
\frac{p(\x)\,\e^{\langle \blambda^*,\h(\x)\rangle }}
{\int_\om \e^{\langle\blambda^*,\h(\x)\rangle}p(\x)\,d\x}\,,\end{equation}
for some vector $\blambda^*=(\lambda^*_j)\in\R^m$; see e.g. \cite[Theorem 3.4]{jordan}. This  vector $\blambda^*$ is an optimal solution of the dual problem described below.

\subsection*{The dual problem}

Associated with the primal convex problem \eqref{def-entropy} is the dual problem
\begin{equation}
 \label{def-entropy-dual}
 \sup_{\blambda\in\R^m} \langle \blambda,\y\rangle -\ln Z(\blambda)\,.
\end{equation}
where one recognizes the Legendre-Fenchel transform $(\ln Z)^*$ of $\ln Z$, evaluated  at $\y$.
Indeed notice that $\blambda^*$ in \eqref{def-optimal-q*} is a critical point of the concave function
$\blambda\mapsto \langle \blambda,\y\rangle -\ln Z(\blambda)$ because
with $q^*$ as in \eqref{def-optimal-q*}:
\begin{eqnarray*}
\y&=&\int_\om \h(\x)\,q^*(\x)\,d\x\\
&=&\frac{\int_\om  \h(\x)\,\e^{\sum_{j=1}^m\lambda^*_j\,h_j(\x)}p(\x)\,d\x}{Z(\blambda^*)}
\,=\,\frac{\nabla Z(\blambda^*)}{Z(\blambda^*)}\,=\,\nabla \ln Z(\blambda^*)\,.
\end{eqnarray*}
Hence when equality between the primal and dual values is obtained,
$\Theta(\y)$ in \eqref{def-entropy} is the Legendre-Fenchel transform $(\ln Z)^*$ in \eqref{def-entropy-dual} 
of the log-partition function $\ln Z(\blambda)$, evaluated at $\y$, i.e.,
\begin{equation}
\label{duality-in-KL}
\Theta(\y)\,=\,
(\ln Z)^*(\y)\,.
\end{equation}
It is also well-known that $\Theta$ is the Cram\'er rate function of the empirical mean of the pushforward 
$\h\#P$ of $P$ in the  LDP literature (see e.g. \cite{dembo,hollander,Yakov}).
For more details on KL-Maxent, its properties and associated computational issues, the interested reader is referred to e.g. 
\cite{bach,lewis,ben-tal,csiszar,dembo,jordan,hollander,Yakov}.

\begin{thm}
\label{th-main}
Let Assumption \ref{ass-1} and Assumption \ref{ass-2} hold, and let $P(d\x)=p(\x)d\x$ be a probability measure 
with support $\overline{\om}$ (hence with density 
$p$ w.r.t. Lebesgue measure). 

If $\y\in\mathrm{int}(\mathcal{M})$ then 
the respective optimal values of \eqref{def-entropy} and the dual \eqref{def-entropy-dual} are attained
at $Q^*(d\x)=q^*d\x$ and $\blambda^*\in\R^m$ with 
$q^*$ as in \eqref{def-optimal-q*}. In addition:

(i) The associated partition function $Z(\blambda^*)$
is the Laplace transform 
evaluated at $-\blambda^*$ 
of the pushforward $\h\#P$ on $\R^m$, with density $v:\R^m\to\R_+$ as in \eqref{def-v}.
That is, $Z(\blambda^*)=\mathcal{L}_v(-\blambda^*)=\mathcal{L}_v^\vee(\blambda^*)$) with $v$ defined in \eqref{def-v}. 

(ii) Hence $\Theta(\y)$  is the \emph{Cram\'er transform} of the function $v$ in \eqref{def-v}, evaluated at $\y$, that is, 
\begin{eqnarray}
\nonumber
\Theta(\y)&=&\sup_{\blambda}\,\langle\blambda,\y\rangle-\ln Z(\blambda)\\
\label{cramer-in-KL}
&=&(\ln Z)^*(\y)\,=\,(\ln \mathcal{L}^\vee_v)^*(\y)
\,=:\,\mathrm{Cramer}_v(\y)\,.
\end{eqnarray}
 \end{thm}
\begin{proof}
The first statement follows from \cite[Theorem 3.4, p. 67]{jordan}.
 Then for (i) and (ii), by Lemma \ref{lem1} with $\blambda=\blambda^*$,
 $Z(\blambda^*)= \mathcal{L}_v(-\blambda^*) =\mathcal{L}^\vee_v(\blambda^*)$
 with $v$ as in \eqref{def-v}.
 Next, by the duality
 \eqref{duality-in-KL} in KL-Maxent, 
 \[\Theta(\y)\,=\,(\ln Z)^*(\y)\,=\,(\ln \mathcal{L}^\vee_v)^*(\y)\,,\]
 which yields the desired result \eqref{cramer-in-KL}.
\end{proof}
So Theorem \ref{th-main} states that the optimal value function $\Theta$  
of the KL-optimization problem \eqref{pb-3} is the Cram\'er Transform of the value 
function $v$ of the integral \eqref{def-v}, a computational geometry problem. This is another instance
of a rigorous link (via the Cram\'er transform) between the well-known KL-Maxent optimization
problem and an integration problem in computational geometry.

If $\y\in\overline{\mathcal{M}}\setminus\mathrm{int}(\mathcal{M})$ then
$\Theta(\y) =\lim_{n\to\infty}\Theta(\y^{(n)})$ with $(\y^{(n)})_{n\in\N}\subset\mathrm{int}(\mathcal{M})$. Finally 
$\Theta(y)=+\infty$ whenever $\y\not\in\overline{\mathcal{M}}$ (see \cite[Theorem 3.4]{jordan}).

One could consider exponential families which are not regular nor minimal,
and to invoke more involved (or abstract) conditions in the literature
(e.g. \cite{lewis,brown,csiszar,gray}). But 
the above framework of regular and minimal exponential families fits
our purpose of illustrating 
how KL-Maxent is related with other area.

\section{Connecting LP and SDP with Max-entropy}

In this section we show that to the canonical linear program (LP) 
one may associate a certain KL-Maxent problem \eqref{pb-3}. Then 
the  perspective function $\varepsilon\,\Theta(\y/\varepsilon)$ associated with the optimal value 
$\Theta(\y)$ of \eqref{pb-3}, converges to the LP optimal value as $\varepsilon\downarrow 0$.
A similar result also holds for the canonical semidefinite program (SDP).

\subsection{Max-entropy and Linear programming  (LP)}

With $\c\in\R^d$, and $\A\in\R^{m\times d}$ with $m<d$, consider the canonical LP 
\begin{equation}
\label{def-lp}
\tau(\y)\,=\,\min_{\x\geq0}\{\,\langle\c,\x\rangle: \A\,\x=\y\,.\}\,\end{equation}
where one may safely assume that $\A$ is full row rank (and so $\A\A^T\succ0$).

Assume that there exists $\blambda_0\in\R^m$ such that $\A^T\blambda_0<\c$. Then
$\tau(\y)\geq\langle \blambda_0,\y\rangle$ for all $\y\in\R^m$. In this case and with no loss of generality,
one may and will assume that $\c>0$; otherwise replace $\c$ with $\tilde{\c}:=\c-\A^T\blambda_0>0$, and $\tilde{\blambda}_0:=0$.
The optimal value
is translated by $-\langle\blambda_0,\y\rangle$ (a constant when $\y$ is fixed).
The LP dual of \eqref{def-lp} reads:
\begin{equation}
\label{def-lp-dual}
\tau(\y)\,=\,\max_{\blambda}\{\,\langle\blambda,\y\rangle: \A^T\blambda\,\leq\,\c\,\}\,,\end{equation}
and in the log-barrier formulation of \eqref{def-lp-dual} for interior point algorithms, one solves
\begin{equation}
\label{lp-barrier}
\tau_\varepsilon(\y)\,=\,\left\{\begin{array}{rl}
\displaystyle\sup_{\blambda}\,\langle\blambda,\y\rangle +\varepsilon\sum_{j=1}^d\ln (\c-\A^T\blambda)_j\,,&\mbox{if $\A^T\blambda<\c$}\\
-\infty\,,&\mbox{otherwise,}\end{array}\right.\,,\end{equation}
for $\varepsilon>0$. When $\varepsilon\downarrow 0$ then $\tau_\varepsilon(\y)\to \tau(\y)$. 
See e.g. \cite{etienne,Monteiro}.

\subsection*{Link with max-entropy}
Let $\om=\mathrm{int}(\R^d_+)\,(=\{\x\in\R^d;\x>0\,\})$, $\h(\x)=\A\,\x$, and as $\c>0$, let $s:=
\prod_{j=1}^d c_j$, and 
\[P(d\x)\,:=\,s\,\e^{\langle -\c,\x\rangle}\,d\x\,,\]
so that $P$ is a probability measure on $\om$
($P(d\x)=\otimes\nu_j(dx_j)$ where $\nu_j$ is the (scaled) exponential
probability measure on $\R_+$, with parameter $c_j$).  Then consider the KL-Maxent problem:
\begin{equation}
\label{def-Theta}
\Theta(\y)\,=\,\inf_{Q\ll P}\,\{\,D_{KL}(Q\Vert P): \int \h\,dQ\,=\,\y\,\}\,.
\end{equation}
To $\Theta$ in \eqref{def-Theta} and $\varepsilon >0$,
is associated the \emph{perspective  function}  \cite{perspective}
\begin{equation}
\label{def-perpsective}
(\y,\varepsilon)\mapsto \widetilde{\Theta}(\y,\varepsilon)\,:=\,
\varepsilon\,\Theta(\y/\varepsilon)\,,
\end{equation}
and let $\widetilde{\Theta}(\y,0):=\lim_{r\downarrow 0}\,r\,\Theta(\y/r)$ 
whenever the limit exists and is finite.
Recall that we assume that $A$ is full row rank. Observe that since
$\overline{\om}=\R^d_+$ and $\h(\x)=\A\x$ then
with $\mathcal{M}$ as in \eqref{moment cone},
$\mathrm{int}(\mathcal{M})=\mathrm{int}(\A\R^d_+)$. Indeed simply use the linearity
$\int_{\R^d_+} \A\x \,\mu(d\x)\,=\,\A\,\int_{\R^d_+} \x\,\mu(d\x) \in\A\R^d_+$, and
if $\mathbf{z}=\A\,\x_0$ for some $\x_0\in\R^d_+$, then $\z=\int_{\R^d_+} \A\x\,\delta_{\x_0}(d\x)\in \mathcal{M}$.
\begin{thm}
\label{th-LP}
Suppose that with $m<d$, $\A^T\blambda_0<\c$ for some $\blambda_0\in\R^m$.
Let $\y\in\mathrm{int}(\A\R^d_+)$.
Then up to the constant $-\varepsilon\,\ln s$, 
$\widetilde{\Theta}(\y,\varepsilon)$
is also the optimal value
$\tau_\varepsilon(\y)$ in \eqref{lp-barrier} of the log-barrier algorithm to solve the LP dual  \eqref{lp-dual}. 
In addition, $\widetilde{\Theta}(\y,0)=\tau(\y)$, the optimal value of the LP
\eqref{def-lp}. 
\end{thm}
\begin{proof}
Observe that $j_\h(\x)\,=\,\sqrt{\mathrm{det}(\A\A^T)}>0$ (independent of $\x$).
Hence if $\A^T\blambda<\c$, then
\begin{eqnarray*}
\frac{s}{\prod_{j=1}^d(\c-\A^T\blambda)_j}&=&s\,\int_{\R^d_+}\e^{\langle -\c+\A^T\blambda,\x\rangle}d\x\,=\,
s\,\int_{\R^d_+}\e^{\langle \A^T\blambda,\x\rangle}\e^{\langle -\c,\x\rangle}d\x\\
&=&\int_\om\e^{\langle\blambda,\h(\x)\rangle} dP(\x)\\
&=&\int_{\R^m}\e^{\langle\blambda,\y\rangle} (\underbrace{\int_{\om\cap\h^{-1}(\y)}\frac{s\,\e^{-\langle\c,\x\rangle}}{\sqrt{\mathrm{det}(\A\A^T)}}d\mathcal{H}^{d-m}(\x)}_{v(\y)})\,d\y\\
&=&\mathcal{L}_v(-\blambda)\,=\,\mathcal{L}^\vee_v(\blambda)\,,
\end{eqnarray*}
where we have applied the coarea formula \eqref{coarea-1}.
Therefore
\begin{equation}
\label{Z(lambda)-polytope}
\ln Z(\blambda)\,=\,\ln \mathcal{L}^\vee_v(\blambda)\,=\,\ln s-\sum_{j=1}^d \ln (\c-\A^T\blambda)_j\,.\end{equation}
Notice that $\mathcal{D}_Z$ in \eqref{domain-D_z} is the open set $\{\blambda: \A^T\blambda<c\}$
and as $\A$ is full rank, 
the exponential family $\exp(\langle \blambda,\h\rangle-\ln Z(\blambda))_{\blambda}$ is regular and minimal.
Moreover, $\mathrm{int}(\mathcal{M})=\mathrm{int}(\A\R^d_+)$. 
Hence  by \eqref{cramer-in-KL}, if $\y\in\mathrm{int}(\A\R^d_+)$,
\begin{eqnarray*}
\Theta(\y)\,=\,(\ln Z)^*(\y)&=&\sup_{\blambda}\,\langle\blambda,\y\rangle-\ln Z(\blambda)\\
&=&-\ln s+\sup\,\{\langle\blambda,\y\rangle+\sum_{j=1}^d\ln (\c-\A^T\blambda)_j: \A^T\blambda<\c\,\}\,,
\end{eqnarray*}
and therefore, with $\varepsilon>0$ one obtains the desired result
\begin{eqnarray*}
\varepsilon\,\Theta(\y/\varepsilon)&=&-\varepsilon\,\ln s+\displaystyle\sup_{\blambda}\,\{\langle\blambda,\y\rangle+\varepsilon\,\sum_{j=1}^d \ln (\c-\A^T\blambda)_j: \A^T\blambda<\c\,\}\\
&=&-\varepsilon\,\ln s+\tau_\varepsilon(\y)\,,\end{eqnarray*}
with $\tau_{\varepsilon}$ as in \eqref{lp-barrier}.  In addition,
\[\widetilde{\Theta}(\y,0)\,=\,\lim_{r\downarrow 0}r\,\Theta(\y/r)\,=\,
\lim_{r\downarrow 0}\tau_r(\y)\,=\,\tau(\y)\,.\]
\end{proof}
Notice that as $j_\h(\x)=\sqrt{\mathrm{det}(\A\A^T)}$ is constant, the function $v$ reads
\begin{equation}
\label{def-vv}
v(\y)\,=\,\frac{s}{\sqrt{\mathrm{det}(\A\A^T)}}
\int_{\{\x\geq0: \A\x=\y\}}\e^{-\langle\c,\x\rangle}\,d\mathcal{H}^{d-m}\,.\end{equation}
In fact, there is an explicit and elegant formula for $v$ in terms of basic ingredients 
(basis $\sigma$ (and its associated matrix $\A_\sigma$), vertex $\x(\sigma)\in\R^d_+$,
and dual vector $\pi_\sigma\in\R^m$) involved in the simplex algorithm to solve the primal LP \eqref{def-lp}. It is due to Brion and Vergne \cite{brion} who use 
the inverse Laplace transform and Cauchy's residue theorem 
to obtain $v(\y)$ in closed form.

Namely with $\y$ fixed, let $\Delta$ be the set of feasible bases $\sigma$, identified with 
an invertible submatrix $\A_\sigma=[\A_{\sigma_1},\ldots, \A_{\sigma_m}]$ of $\A$, 
and a vertex $\x(\sigma)\in\R^d$ with $\x(\sigma)_j=(\A_\sigma^{-1}\y)_j$, $j\in\sigma$, and $0$ 
if $j\not\in\sigma$. Introduce
the vector $\c_\sigma=(c_{\sigma_1},\ldots,c_{\sigma_m})\in\R^m$, 
and the dual vector $\pi_\sigma\in\R^m$ defined by $\A^T_\sigma\pi_\sigma=\c_\sigma$. 
Then by \cite[Theorem, p. 801]{brion}, when $\c$ is a generic cost vector and
$\y$ is such that the polyhedron  $\{\x\geq0: \A\x=\y\}$ is simple\footnote{With $\A\in\R^{m\times d}$, $\y\in\R^m$, the polytope $\{\,\x\in\R^d_+:\A\x\,=\,\y\,\}$ is simple if at any of its vertices, exactly $d-m$ coordinates vanish.}, then:
\begin{eqnarray}
\nonumber
v(\y)&=&s\,\sum_{\sigma\in\Delta}\frac{\e^{-\langle \c,\x(\sigma)\rangle}}
{\vert\mathrm{det}(\A_\sigma)\vert\prod_{j\not\in \sigma}(c_j-\pi_\sigma^T\A_j)}\\
\label{brion&vergne}
&=&s\,
\sum_{\sigma\in\Delta}\frac{\e^{-\langle \pi_\sigma,\y\rangle}}
{\vert\mathrm{det}(\A_\sigma)\vert\prod_{j\not\in \sigma}(c_j-\pi_\sigma^T\A_j)}\,.
\end{eqnarray}
The above formula $(4.8)$ is \cite[Theorem p. 801(ii)]{brion} in which
$h\in \bar{\gamma}$ (where $\gamma$ is a \emph{chamber}) 
translates to $\y\in\bar{\gamma}$. 
The set $\mathcal{B}(\Delta,\gamma)$ is the set of all bases $\sigma$ such that
(i) $\A_\sigma=(\A_j)_{j\in\sigma}$ is nonsingular, and (ii)  $\gamma$ is an open set contained  in the cone generated by the columns $(\A_j)_{j\in\sigma}\subset \mathbb{R}^m$. As $\gamma$ is open,  $\A_\sigma^{-1}\y>0$ whenever $\y\in\gamma$. 
In \eqref{brion&vergne} the triplet $(\y, \x(\sigma),\c)$ corresponds to the triplet $(h,v_\sigma(h),y)$ in \cite[Theorem (ii), p. 801]{brion}. In a chamber, all the associated polytopes have the same combinatorial structure. So whenever the polytope $\{\x\geq0:\A\x=\y\}$ is simple, $\y$ belongs to
a unique chamber $\gamma$;
then $\mathcal{B}(\Delta,\gamma)$ in \cite[Theorem (ii), p. 801]{brion} is $\Delta$ in
\eqref{brion&vergne}.

When the polytope is not simple, some vertices are degenerate and $\y$ lies on the boundary of one or several chambers. One then chooses a chamber $\gamma$ such that $\y\in\overline{\gamma}$, and uses the corresponding set $\mathcal{B}(\Delta,\gamma)$ in \cite[Theorem (ii), p. 801]{brion}.
\begin{ex}
Consider the following illustrative toy example with $0<\c=(c_1,c_2)$ and
$(d,m)=(2,1)$ with $\A=[1,1]$. The natural volume element on $\h^{-1}(y)$ is the 
``length" on the line $\{\x:  x_1+x_2=y\}$ of $\R^2$. 
That is, if one parametrizes 
$\R^d_+\cap \h^{-1}(y)$ by $(t,y-t)$, $t\in [0,y]$,  then 
$\mathcal{H}^{d-m}(d\x)=\sqrt{2}dt$ and $v$ in \eqref{def-vv} reads
\[v(y)\,=\,\frac{1}{\sqrt{\mathrm{det}(\A \A^T)}}\int_{\{\x\geq0: x_1+x_2=y\}}s\,\e^{\langle-\c,\x\rangle}d\mathcal{H}^{d-m}\,,\quad y\geq0\,,\]
\[=\,
\frac{c_1c_2}{\sqrt{2}}\int_0^y\e^{-c_1\,y}\e^{(c_1-c_2)t}\sqrt{2}dt\,=\,
\frac{c_1c_2}{c_1-c_2}\,(\e^{-c_2y}-\e^{-c_1y})\,,\]
for all $y\geq0$.
\end{ex}

\subsection{Link with Semidefinite Programming (SDP)}

Let $\mathcal{S}^d$ be the vector space of $(d\times d)$ real symmetric matrices (hence of dimension
$s_d:=d(d+1)/2$) (with usual scalar product $\langle \X,\Y\rangle=\mathrm{trace}(\X\,\Y)$) and let $K:=\mathcal{S}^d_+$ be the convex cone of real $(d\times d)$ positive semidefinite matrices.
Then the dual cone $K^*$ of $K$ coincides with $K$ ($K^*=K$) and its
associated usual log-barrier  function $\varphi:\mathrm{int}(K^*)\to \R$, 
reads
\begin{equation}
\label{logbarrier-sdp}
\varphi(\V)\,:=\,\ln\left(\int_{K}\e^{-\langle \V,\X\rangle}d\X\right)\,=\,C_d-\frac{d+1}{2}\ln \mathrm{det}(\V)\,,\end{equation}
for all $\V\in\mathrm{int}(K)$, and
for some constant $C_d$ that depends only on $d$. (Up to a scaling,
the density $\exp(-\langle \V,\X\rangle) d\X$ is related to the Wishart distribution 
with parameter $(\V^{-1}/2,(d+1)$) on $K$.)

\subsection*{A canonical SDP}
With $m<s_d$, let $\A_0,\ldots,\A_m\in \mathcal{S}^d$ be given, and consider the semidefinite program (SDP):
\begin{equation}
\label{sdp-primal}
\tau(\y)\,=\,\inf_{\X\in K} \{\,\langle \A_0,\,\X\rangle:  \langle \A_j,\X\rangle\,=\,y_j\,,\:j=1,\ldots,m\,\}\,,\end{equation}
with associated dual
\begin{equation}
\label{sdp-dual}
\tau^*(\y)\,=\,\sup_{\blambda\in\R^m} \{\,\langle \blambda,\y\rangle:  \A_0-\sum_{j=1}^m\lambda_j\,\A_j\in K^*\,\}.\end{equation}
With no loss of generality we may and will assume that $\A_1,\ldots,\A_m$
are linearly independent in $\mathcal{S}^d$. We also assume that 
\begin{equation}
\label{slater-sdp}
\exists \blambda_0\in\R^m:\quad \A_0-\sum_{j=1}^m\lambda_{0j}\,\A_j\,\succ 0\,,\end{equation}
that is, Slater's condition holds for \eqref{sdp-dual} and therefore $\tau(\y)=\tau^*(\y)$
and \eqref{sdp-primal} has an optimal solution $\X^*\in K$ whenever $\tau(\y)$ is finite. 

Hence under \eqref{slater-sdp} and with no loss of generality,
we may and will assume that $\A_0\succ0$; otherwise replace $\A_0$ with $\tilde{\A}_0:=
\A_0-\sum_{j=1}^m\lambda_{0j}\,\A_j\,\succ 0$, so that Slater's condition is satisfied with $\tilde{\blambda}_0:=0$, and the optimal value is translated by the constant $-\langle\blambda_0,\y\rangle$.

A  classical approach to solve \eqref{sdp-dual} is via its associated log-barrier formulation
with $\varepsilon$-parameter, for interior point algorithms.
Namely, with $\varepsilon>0$, solve
\begin{equation}
\label{sdp-dual-barrier}
\tau^*_\varepsilon(\y)\,=\,\sup_{\blambda\in\R^m} \{\,\langle \blambda,\y\rangle -\varepsilon\,\varphi(\A_0-\sum_{j=1}^m\lambda_j\,\A_j)
\end{equation}
where the sup is over all $\blambda$ such that $\A_0-\sum_{j=1}^m\lambda_j\,\A_j\in K$. Then
$\tau^*_\varepsilon(\y)\to \tau^*(\y)$ as $\varepsilon\downarrow 0$.
See e.g. \cite{etienne,log-barrier-3}.

\subsection*{Link with max-entropy}
Let $\om:={\color{blue}\mathrm{int}(K)}$, and since $\A_0\succ0$, let $1/s:=\int_K\e^{-\langle \A_0,\X\rangle}d\X$, and introduce
the probability measure $P(d\X)$ on $\om$,
\begin{equation}
\label{sdp-p}
P(d\X)\,=\,p(\X)\,d\X\,=\,s\,\e^{-\langle \A_0,\X\rangle}d\X\,.
\end{equation}
In this context, let $\h:\mathcal{S}^d\to\R^m$
be the linear mapping 
\begin{equation}
\label{sdp-h}
\X\mapsto  \h(\X)\,:=\,(\langle \A_j,\X\rangle)_{j=1,\ldots,m}\,,\quad\forall \X\in \mathcal{S}^d\,,\end{equation}
and with  $J_\h(\X)$ being the Jacobian of $\h$,
$J_\h(\X)J_\h(\X)^T$ is a $(m\times m)$-matrix with entries 
\[J_\h(\X)J_\h(\X)^T(i,j)\,=\,\mathrm{trace}(\A_i\cdot \A_j^T)\,,\quad 1\leq i,j\leq m\,.\]
Therefore
\begin{equation}
\label{sdp-jh}
j_\h(\X)\,=\,\sqrt{\mathrm{det}(J_\h(\X)J_\h(\X)^T)}\,=\,\mbox{Cte}=:\,C_{\A_1,\ldots,\A_m}\,.
\end{equation}
Moreover $C_{\A_1,\ldots,\A_m}>0$ because the $\A_j$'s are linearly independent in $\mathcal{S}^d$.
Indeed let $\mathbf{G}:=(\langle \A_i,\A_j\rangle)_{1\leq i,j\leq m}$,
and $\z\in \R^d$, $\z\neq0$. Then:
\begin{eqnarray*}
\langle \z,\mathbf{G}\z\rangle&=&\sum_{i,j}z_i\,z_j\,\langle \A_i,\A_j\rangle\\
&=&\langle\sum_{i=1}^m\A_i\,z_i,\sum_{j=1}^m\A_j\,z_j\rangle\,=\,
\Vert\sum_{i}\A_i\,z_i\Vert^2\,>0\,,
\end{eqnarray*}
where the last inequality follows by linear independence of the $\A_j$.

Next, consider the KL-Maxent problem
\begin{equation}
\label{entropy-sdp}
\begin{array}{rl}
\Theta(\y)=\displaystyle\min_{Q\ll P}&\{\,\displaystyle\int_\om q(\X)\,\ln \frac{q(\X)}{p(\X)}d\X:\\
\mbox{s.t.}&  \displaystyle\int_\om \langle \A_j,\X\rangle\,dQ(\X)\,=\,y_j\,,\quad j=1,\ldots,m\,\}\,
\end{array}
\end{equation}
(with $Q(\X)=q(\X)d\X$). Define $v:\R^m\to \R$ by:
\begin{equation}
\label{v-sdp}
\y\mapsto v(\y)\,:=\,\frac{s}{C_{\A_1,\ldots,\A_m}}\,\int_{K\cap \h^{-1}(\y)}\e^{\langle -\A_0,\X\rangle}d\X\,,\end{equation}
with $v(\y)=0$ if $\om\cap\h^{-1}(\y)=\emptyset$.
\begin{lem}
\label{lem-sdp}
Assume that $j_\h>0$. Then for real $\blambda\in\R^m$ such that $\A_0+\sum_{j=1}^m\lambda_j\,\A_j\in \mathrm{int}(K)$,
the Laplace transform $\mathcal{L}_v(\blambda)$ of $v$ reads:
\begin{equation}
\mathcal{L}_v(\blambda)
\label{lem-sdp-1}
\,=\,s\,\e^{\varphi(\A_0+\sum_{j=1}^m \lambda_j\,\A_j)}\,,
\end{equation}
where $\varphi$ is the log-barrier in \eqref{logbarrier-sdp} associated with  $K^*=K$. 
\end{lem}
\begin{proof}
As $\A_0+\sum_{j=1}^m\lambda_j\,\A_j\in \mathrm{int}(K)$,
 \begin{eqnarray*}
 \mathcal{L}_v(\blambda)&=&\int_{\R^m}\e^{-\langle\blambda,\y\rangle}v(\y)\,d\y\\
 &=&\int_{\R^m}\e^{-\langle\blambda,\y\rangle}\left(\frac{s}{C_{\A_1,\ldots,\A_m}}\,\int_{K\cap\h^{-1}(\y)}\e^{-\langle \A_0,\X\rangle}d\X\right)d\y\\
 &=&\int_K\e^{-\langle\blambda,\h(\X)\rangle}\,dP(\X)\quad\mbox{[by \eqref{coarea-1}]}\\
 &=&s\,\int_K\e^{-\langle\blambda,\h(\X)\rangle}\,\e^{-\langle \A_0,\X\rangle}d\X\,=\,s\,\e^{\varphi(\A_0+\sum_{j=1}^m\lambda_j\,\A_j)}\,.
 \end{eqnarray*}
 \end{proof}
 Observe that with $\mathcal{M}$ as in \eqref{moment cone}, $\h(K)=\mathcal{M}$. Indeed
 for any $\mu\in \mathscr{P}(K)$, and by linearity of $\h$,
  $\int_K \h(\X)\,\mu(d\X)=\h(\int_K \X\,\mu(d\X))\in\h(K)$. Conversely, if
  $\Z=\h(\X_0)$ for some $\X_0\in K$, then $\Z=\int_K \h(\X)\,\delta_{\X_0}(d\X)\in\mathcal{M}$.
\begin{thm}
\label{thm-sdp}
Let $\om:={\color{blue}\mathrm{int}(K)}$, $\h$ be as in \eqref{sdp-h} with $m\leq d$, and with $\A_0\in \mathrm{int}(K)$,
let $P(d\X)$ be the probability measure on $\om$ defined in \eqref{sdp-p}. Assume that 
$j_\h=C_{\A_1,\ldots,\A_m}$ in \eqref{sdp-jh} is strictly positive, and let 
$\y\in \mathrm{int}(\h(K))$ ($=\mathrm{int}(\mathcal{M})$).

With $\Theta(\y)$ being the optimal value of the Max-Ent problem \eqref{entropy-sdp}, its associated 
perspective function $(\y,\varepsilon)\mapsto \widetilde{\Theta}(\y,\varepsilon)$ ($\varepsilon>0$) is also
(up to the constant $-\varepsilon\,\ln s$) the optimal value of \eqref{sdp-dual-barrier} 
in the log-barrier interior point algorithm associated with the semidefinite program \eqref{sdp-dual}.
In particular, $\tau^*(\y)\,=\,\widetilde{\Theta}(\y,0)\,(:=\lim_{r\downarrow 0}r\,\widetilde{\Theta}(\y/r)$.
\end{thm}
\begin{proof}
 From the results in \S \ref{sec:main} and Lemma \ref{lem-sdp}, 
\begin{eqnarray*}
\ln Z(\blambda)&=&\ln \int_\om e^{\langle\blambda,\h(\X)\rangle}\,dP(\X)\\
&=&\ln \mathcal{L}_v(-\blambda)\,=\,
\ln s+\varphi(\A_0-\sum_{j=1}^m\lambda_j\,\A_j)\,.
\end{eqnarray*}
Next, as $\y\in\mathrm{int}(\mathcal{M})$, then by \eqref{cramer-in-KL}:
\begin{eqnarray*}
\Theta(\y)
&=&\displaystyle\sup_{\blambda\in\R^m}\,\langle \blambda,\y\rangle-\ln Z(\blambda)\\
&=&\displaystyle\sup_{\blambda\in\R^m}\,\langle \blambda,\y\rangle-\ln\mathcal{L}_v(-\blambda)\\
&=&-\ln s+\displaystyle\sup_{\blambda\in\R^m}\,\langle \blambda,\y\rangle-\phi(\A_0-\sum_{j=1}^m\lambda_j\,\A_j)\,.
\end{eqnarray*}
Hence for every $\varepsilon>0$ fixed, arbitrary, 
\begin{equation}
\widetilde{\Theta}(\y,\varepsilon)\,=\,\varepsilon\,\Theta(\y/\varepsilon)\,=\,-\varepsilon\,\ln s+\displaystyle\sup_{\blambda}\,\langle\blambda,\y\rangle-\varepsilon\,\varphi(\A_0-\sum_{j=1}^m\lambda_j\,\A_j)\,,
\end{equation}
where the sup is over all $\blambda$ such that $\A_0-\sum_{j=1}^m\lambda_j\,\A_j\in K$,
and $\widetilde{\Theta}(\y,\varepsilon)$ (with $\varepsilon>0$), is the perspective function associated with $\Theta$.
In particular $\tau^*(\y)\,=\,\widetilde{\Theta}(\y,0)\,(:=\lim_{r\downarrow 0}r\,\tilde{\Theta}(\y/r))$.
\end{proof}

\section{Conclusion}

In this paper we have established a rigorous link between 
on the one hand, the max-entropy optimization problem with KL-divergence,
and on the other hand, an explicit integral in computational geometry, namely a certain integration of 
the density of the reference distribution. The link is obtained by 
geometric measure arguments and in doing so, the optimal value function 
in KL-optimization is nothing less than the Cram\'er transform of the latter integral.
In addition, one has also shown that the optimal value function of the canonical LP
is the limit of the perspective function of a related max-entropy problem with appropriate
reference distribution, and similarly for the canonical SDP.

Links between optimization and integration have already been revealed and discussed in several contexts, as formally the same operation in 
different algebra. In the present context of max-entropy optimization,
its formal link with integration is rigorously established via the Cram\'er transform
of a certain explicit integral in computational geometry.


\begin{thebibliography}{lass}
\bibitem{linearity}
F. Bacelli, G. Cohen, G.J. Olsder, and J.P. Quadrat, \emph{Synchronization and Linearity: An Algebra for Discrete Event Systems}, John Wiley \& Sons Ltd., 1992
\bibitem{bach}
F. Bach. Sum-of-Squares relaxations for information theory and variational inference,
\emph{Found. Comp. Math.} {\bf 25}, pp. 865--903, 2024.
\bibitem{BC}
H.H. Bauschke, P.L. Combettes. \emph{Convex Analysis and Monotone Operator Theory in Hilbert Spaces}, Springer, 2nd ed., 2019
\bibitem{ben-tal}
A. Ben-Tal, M. Teboulle, A. Charnes. The role of duality in optimization problems involving Entropy functionals with applications to Information Theory, \emph{J. Optim. Theory \& Appl.} {\bf 58}(2), pp. 209--223, 1988.
\bibitem{lewis}
J. Borwein, A. Lewis. Duality relationship for entropy-like minimization problems, \emph{SIAM J. Control Optim.} {\bf 29}(2), pp. 325--338, 1991.
\bibitem{brion}
M. Brion, N. Vergne. Residue formulae, vector partition functions, and  lattice points in rational polytopes,
\emph{J. Amer. Math. Soc.} {\bf 10}(4), pp. 797--833, 1997.
\bibitem{brown}
L.D. Brown. \emph{Fundamentals of statistical exponential families: {\small with applications in statistical decision theory}}, Lecture Notes-Monograph Series, Institute of Mathematical Statistics, vol 9, 1986.
\bibitem{perspective}
P.L. Combettes, C.L. M\"uller. Perspective functions: Proximal calculus and applications
in high-dimensional statistics, \emph{J. Math. Anal. Appl.} {\bf  457}, pp. 1283--1306, 2018
\bibitem{Cramer}
H. Cram\'er. Sur un nouveau th\'eor\`eme-limite de la th\'eorie des probabilit\'es,
\emph{Actual. Sci. Ind.} {\bf 736}, pp. 5--23, 1938.
\bibitem{csiszar}
I. Csisz\'ar. I-Divergence geometry of probability distributions and minimization problems, \emph{Ann. Prob.} {\bf 3}(1), pp. 146--158,1975.
\bibitem{dembo}
A. Dembo, O. Zeitouni. \emph{Large Deviations Techniques and Applications}, Springer-Verlag, Berlin, 2009.
\bibitem{Donsker}
M.D. Donsker, S.R.S. Varadhan. Asymptotic evaluation of certain Markov process expectations for large time --III,
\emph{Comm. Pure Appl. Math.} {\bf 29}(4), pp. 389--461, 1976. 
\bibitem{federer}
H. Federer. Curvature measures, \emph{Trans. Amer. Math. Soc.} 
{\bf 93}(3), pp. 418--491, 1959.
\bibitem{federer-book}
H. Federer. \emph{Geometric Measure Theory}, Springer series Classics in Mathematics, 
Springer Berlin, Heidelberg, 1996.
\bibitem{paper-4}
F. Gamboa, C. Gueneau, T. Klein, E. Lawrence.
Maximum Entropy on the mean approach to solve generalized inverse problems with an application in computational thermodynamics,
\emph{RAIRO-Oper. Res.} {\bf 55}, pp. 355--393, 2021.
\bibitem{log-barrier-3}
L.M. Gra\~{n}a Drummond, A.N. Iusem, B.F. Svaiter. On the central path for nonlinear semidefinite
programming, \emph{RAIRO. Oper. Res.} {\bf 34}(3), pp. 331--345, 2000.
\bibitem{gray}
R.M. Gray. \emph{Entropy and Information Theory}, 2nd ed., Springer, New York, 2011.
\bibitem{etienne}
M. Halick\'a, E. de Klerk, C. Roos.  On the convergence of the central path in semidefinite optimization,
\emph{SIAM J. Optim.} {\bf 12}(4), pp. 1090--1099, 2002
\bibitem{didier}
D. Henrion. Maximal entropy in the moment body, {\tt arXiv:2507.02461}, 2025.
\bibitem{hollander}
F. den Hollander.  \emph{Large Deviations}, Fields Institute Monographs {\bf 14}, American Mathematical Society, 2000.
\bibitem{paper-7}
M. King-Roskamp, R. Choksi, T. Hoheisel.
Data-Driven Priors in the Maximum Entropy on the Mean Method for Linear Inverse Problems. {\tt arXiv:2412.17916}, 2024.
\bibitem{lasserre-book}
J.B. Lasserre. \emph{Linear an Integer Programming vs Linear Integration and Counting},
Springer-Verlag, New York, 2009.
\bibitem{Monteiro}
R.D.C. Monteiro, F. Zhou. On the existence and convergence of the central path for convex programming and some duality results,
\emph{Comput. Optim. Appl.} {\bf 10}(1), pp. 51--77, 1998
\bibitem{paper-5}
J. Navaza. On the Maximum-Entropy estimate of the electron density function,
\emph{Acta Cryst.} {\bf  A41}, pp. 232--244, 1985.
\bibitem{paper-6}
J. Navaza. The use of non-local constraints in Maximum-Entropy
electron density reconstruction, \emph{Acta Cryst.} {\bf  A42}, pp. 212--223, 1986.
\bibitem{Negro}
L. Negro. Sample distribution theory using Coarea Formula,
\emph{Communication in Statistics: Theory and Methods} {\bf 53}(5), pp. 1864--1889, 2024.
\bibitem{paper-2}
G. Rioux, R. Choksi, T. Hoheisel, P. Mar\'echal, C. Scarvelis.
The Maximum Entropy on the Mean method for image deblurring, \emph{Inverse Problems} {\bf 37}(1), 015011, 2021.
\bibitem{paper-3}
G. Rioux, C. Scarvelis, R. Choksi, T. Hoheisel, P. Mar\'echal.
Blind deblurring of barcodes via Kullback-Leibler Divergence,
\emph{IEEE Trans. Pattern Anal. Machine Intelligence} {\bf 43}(1), pp. 77--88, 2021
 \bibitem{Schmudgen}
K. Schm\"udgen. \emph{The Moment Problem}, Graduate Texts in Mathematics, Springer Cham, 2017.
 \bibitem{Straszak}
 D. Straszak, N.K. Vishnoi. Maximum Entropy Distributions: Bit Complexity and Stability, Proc. Mach. Learning Research {\bf 99}, pp. 1--31, 2019  
 \bibitem{jordan}
M.J. Wainwright, M.I. Jordan.  
Graphical models, exponential families, and variational inference,
\emph{Found. Trends Mach. Learn} {\bf 1}(1-2), pp. 1--305, 2008
\bibitem{Yakov}
Y. Vaisbourd, R. Choksi, A. Goodwin, T. Hoheisel, C.-R. Sch\"onlieb. Maximum Entropy on the mean and the Cram\'er rate function in statistical estimation and inverse problems: Properties, models and algorithms, \emph{Math. Program.} 
{\bf 214}, pp; 441--490, 2025.
\end{thebibliography}
\end{document}